\documentclass[11pt]{article}
\usepackage{amssymb,amsfonts,amsmath,latexsym, epsfig,mathrsfs,colordvi, xcolor,tikz}
\usepackage{ytableau}
\usepackage{graphicx}
\newtheorem{theo}{Theorem}

\newtheorem{exam}[theo]{Example}
\newtheorem{lemma} [theo]{Lemma}
\newtheorem{coro}[theo]{Corollary}
\newtheorem{prop}[theo]{Proposition}

\makeatletter \@addtoreset{equation}{section}
\@addtoreset{theo}{}\makeatother

\def\qed{\hfill \rule{4pt}{7pt}}
\def\pf{\noindent {\it Proof.} }

\def\Inc{{\mathrm{Inc}}}
\def\RInc{{\mathrm{RInc}}}
\def\pRInc{{\mathrm{pRInc}}}
\def\SYT{{\mathrm{SYT}}}
\def\maj{{\mathrm{maj}}}
\def\amaj{{\mathrm{amaj}}}

\begin{document}

\title{Enumeration on row-increasing tableaux of shape $2 \times n$}

\author{Rosena R. X. Du\footnote{Corresponding Author. Email: rxdu@math.ecnu.edu.cn.}, Xiaojie Fan and Yue Zhao\\ \\ School of Mathematical Sciences, Shanghai Key Laboratory of PMMP \\East China Normal University,
500 Dongchuan Road \\Shanghai, 200241, P. R. China.}

\date{December 21,  2018}
\maketitle
\noindent {\bf Abstract:}
Recently O. Pechenik studied the cyclic sieving of increasing tableaux of shape $2\times n$, and obtained a polynomial on the major index of these tableaux, which is a $q$-analogue of refined small Schr\"{o}der numbers. We define row-increasing tableaux and study the major index and amajor index of row-increasing tableaux of shape $2 \times n$. The resulting polynomials are both $q$-analogues of refined large Schr\"{o}der numbers. For both results we give bijective proofs.



\section{Introduction}

Let $n$ be a positive integer and $\lambda$ a partition of $n$. A \emph{semistandard (Young) tableau} (SSYT) of shape $\lambda$ is an array $T$ of positive integers of shape $\lambda$ that is strictly increasing in every row and weakly increasing in every column. If an SSYT of shape $\lambda$ is strictly increasing in each row and column, and the entries are $1, 2, \ldots, n$, it is called a \emph{standard Young tableau} (SYT). (Note that in many other literatures such as \cite{EC2}, an SSYT is defined to be strictly increasing in every column and weakly increasing in every row instead, but we prefer to define it in this way for later convenience.) We denote by $\SYT(\lambda)$ the set of standard Young tableaux of shape $\lambda$. Throughout this paper we will identify a partition $\lambda$ with its Young diagram; hence the notations $\SYT(m \times n)$ and $\SYT(n^m)$ are equivalent.

A \emph{descent} of an SSYT $T$ is 
any instance of $i$ followed by an $i+1$ in a lower row of $T$, and define the \emph{descent set} $D(T)$ to be the set of all descents of $T$. The \emph{major index} of $T$ is defined by $\maj(T)=\sum_{i \in D(T)}i$. Similarly, we define an \emph{ascent} of $T$ to be any instance of $i$ followed by an $i+1$ in a higher row of $T$, and define the \emph{ascent set} $A(T)$ to be the set of all ascents of $T$. The \emph{amajor index} of $T$ is defined by $\amaj(T)=\sum_{i \in A(T)}i$. For example, for the tableau $T$ shown in Figure \ref{fig:Inc_SYT}, we have $D(T)=\{2,5,6,8\}$, $\maj(T)=21$, $A(T)=\{3,4,7\}$ and $\amaj(T)=14$.

The most well-known result on enumerating the major indices of Young tableaux is the following \emph{$q$-hook length formula}:

\begin{lemma} \cite[p.376]{EC2} For any partition $\lambda=\sum_{i}\lambda_{i}$ of $n$, we have
\begin{equation}\label{qhooklength}
\sum_{T \in \SYT(\lambda)}q^{\maj(T)}=\frac{q^{b(\lambda)}[n]!}{\prod_{u \in \lambda}h(u)}.
\end{equation}
Here $b(\lambda)=\sum_{i}(i-1)\lambda_{i}$.
\end{lemma}

Applying the above result to SYTs of shape $2 \times n$, we get the following result.
\begin{equation}\label{Cqn}
C_q(n)=\sum_{T \in \SYT(2 \times n)}q^{\maj(T)}=\frac{q^n}{[n+1]}{2n \brack n}.
\end{equation}
The above result is a $q$-analogue of the well-known result that SYTs of shape $2 \times n$ are counted by the $n$-th \emph{Catalan number} $C_n=\frac{1}{n+1} {2n \choose n}$.
The famous RSK algorithm \cite{EC2} is a bijection between permutations of length $n$ and pairs of SYT of order $n$ of the same shape. Under this bijection, the descent set of a permutation is transferred to the descent set of the corresponding ``recording tableau”. Therefore many problems involving the statistic \emph{descent} or \emph{major index} of pattern-avoiding permutations can be translated to the study of descent and major index of standard Young tableaux \cite{BarnabeiBonettiElizaldeSilimbani,KeithMajor}.

In this paper we study the major (amajor) index polynomial of increasing and row-increasing tableaux. An \emph{increasing tableau} is an SSYT such that both rows and columns are strictly increasing, and the set of entries is an initial segment of positive integers (if an integer $i$ appears, positive integers less than $i$ all appear).  We denote by $\Inc_{k}(\lambda)$ the set of increasing tableaux of shape $\lambda$ whose entries are $\{1, 2, \ldots, n-k\}$, i.e., $\Inc_{k}(\lambda)$ denotes the set of increasing tableaux of shape $\lambda$, with exactly $k$ numbers appear twice.  Figure \ref{fig:Inc_SYT} shows an increasing tableau $T \in \Inc_3( 2 \times 6)$.

\begin{figure}[ht]
\centering
$T:$\ \ \ytableaushort{1{\color{blue} 2} {\bf 4}{\color{blue} 5}{\color{blue}\bf 6}{\color{blue}\bf 8},
{\color{red} 3} {\color{red}\bf 4} {\bf 6} {\color{red} 7} {\bf 8}9}
\caption{An increasing tableau $T \in \Inc_3( 2 \times 6)$.}\label{fig:Inc_SYT}
\end{figure}

In \cite{OliverP_JCTA} O. Pechenik studied increasing tableaux in $\Inc_k(2 \times n)$ and obtained the following result.
\begin{theo}\cite[Pechenik]{OliverP_JCTA}\label{theo_sqnk}  For any positive integer $n$, and $0 \leq k \leq n$ we have
\begin{equation}\label{sqnk}
S_q(n,k)=\sum_{T \in \Inc_k(2 \times n)}q^{\maj(T)}=\frac{q^{n+k(k+1)/2}}{[n+1]}{n-1 \brack k}{2n-k \brack n}.
\end{equation}
\end{theo}

Note that $\Inc_0(2 \times n)=\SYT(2 \times n)$, and Equation \eqref{sqnk} coincides with Equation \eqref{Cqn} when $k=0$. Moreover, setting $q=1$, Equeation \eqref{sqnk} indicates that the cardinality of $\Inc_k(2 \times n)$ is
\begin{equation}\label{snk}
s(n,k)=\frac{1}{n+1}{n-1 \choose k}{2n-k \choose n}.
\end{equation}
$s(n,k)$ is sequence A126216 in OEIS \cite{sequence}, which is considered as a refinement of the \emph{small Schr\"{o}der number} (sequence A001003 in OEIS \cite{sequence}). And $s(n,k)$ also counts the following sets.
\begin{itemize}
\item[1.] Dissections of a convex $(n+2)$-gon into $n-k$ regions;
\item[2.] SYTs of shape $(n-k,n-k,1^k)$;
\item[3.] Noncrossing partitions
of $2n-k$ into $n-k$ blocks all of size at least 2.
\end{itemize}

In \cite{StanPolySYT} Stanley gave a bijection between 1 and 2. In \cite{OliverP_JCTA} Pechenik gave bijections between 2, 3 and increasing tableaux in $\Inc_k(2 \times n)$. Moreover,  Pechenik’s bijection between $\Inc_k(2 \times n)$ and $\SYT(n-k, n-k, 1^k)$ preserves the descent set, therefore by applying the $q$-hook length formula to $\SYT(n-k, n-k, 1^k)$ Pechenik proved Theorem \ref{theo_sqnk}. 

In the theory of lattice paths enumeration $s(n,k)$ also counts the number of small Schr\"{o}der $n$-paths with $k$ $F$ (flat) steps. Here a \emph{Schr\"{o}der $n$-path} is a lattice path goes from $(0,0)$ to $(n,n)$ with steps $(0,1)$($U$), $(1,0)$($D$) and $(1,1)$($F$) and never goes below the diagonal line $y=x$. If a  Schr\"{o}der path contains no $F$ steps on the diagonal line, it is called a \emph{small Schr\"{o}der path}. There is an obvious bijection between SSYTs in $\Inc_k(2 \times n)$ and small Schr\"{o}der $n$-paths with $k$ flat steps: read the numbers in a tableau from $1$ to $2n-k$ in increasing order, if $i$ appears only in row 1 (2), it corresponds to a $U$ ($D$) step; if $i$ appears in both rows, it corresponds to an $F$ step. 
Increasing rectangular tableaux and $m$-Schr\"{o}der paths are also studied in \cite{PresseyStokkeVisentin, SongCW_GernalizedSchroder}.

Note that the above bijection works for all Schr\"{o}der $n$-paths, and an $F$ step on the diagonal line corresponds to a column with identical numbers in a tableau of shape $2 \times n$. Motivated by Pechenik’s results and the above observation, we define and study \emph{row-increasing tableaux}. Here a row-increasing tableau is an SSYT with strictly increasing rows and weakly increasing columns, and the set of entries is a consecutive segment of positive integers.  Given positive integer $n$, nonnegative integers $k, m$, and $\lambda$ a partition of $n$, we denote by $\RInc^m_{k}(\lambda)$ the set of row-increasing tableaux of shape $\lambda$ with set of entries $\{m+1, m+2, \ldots, m+n-k\}$. We will also denote $\RInc^0_{k}(\lambda)$ as $\RInc_{k}(\lambda)$. It is obvious that $\Inc_{k}(\lambda) \subseteq \RInc_{k}(\lambda)$.

It is not hard to see (will be explained in Section \ref{sec2.3}) that $\RInc_{k}(2 \times n)$ is in bijection with Schr\"{o}der $n$-paths with $k$ $F$ steps, and these two sets are both counted by
\begin{equation}\label{rnk}
r(n,k)=\frac{1}{n-k+1}{2n-k \choose k}{2n-2k \choose n-k}.
\end{equation}
Here $r(n,k)$ is considered as a refinement of the \emph{large Schr\"{o}der number} and is the sequence A006318 in OEIS \cite{sequence}.

Our main results are the following formulas involving major index and amajor index of SSYTs in $\RInc_{k}(2 \times n)$.
\begin{theo}\label{theo_rqnk} For any positive integer $n$, and $0 \leq k \leq n$ we have
\begin{equation}\label{rqnk}
R_q(n,k)=\sum_{T \in \RInc_k(2 \times n)}q^{\maj(T)}=\frac{q^{n+k(k-3)/2}}{[n-k+1]}{2n-k \brack k}{2n-2k \brack n-k}.
\end{equation}
\end{theo}

\begin{theo}\label{theo_r2qnk} For any positive integer $n$, and $0 \leq k \leq n$ we have
\begin{equation}\label{r2qnk}
\widetilde{R}_q(n,k)=\sum_{T \in \RInc_k(2 \times n)}q^{\amaj(T)}=\frac{q^{k(k-1)/2}}{[n-k+1]}{2n-k \brack k}{2n-2k \brack n-k}.
\end{equation}
\end{theo}

Note that when $k=0$, we have $\RInc_0(2 \times n)=\SYT(2 \times n)$. Therefore Theorem \ref{theo_rqnk} indicates Equation \eqref{Cqn} and Theorem \ref{theo_r2qnk} indicates the following result.
\begin{coro} For any positive integer $n$, we have
\begin{equation}
\widetilde{C}_q(n)=\sum_{T \in \SYT(2 \times n)}q^{\amaj(T)}=\frac{1}{[n+1]}{2n \brack n}.\label{c2qnk}
\end{equation}
\end{coro}

The organization of the paper is as follows. In Section 2 we give a bijection between row-increasing tableaux in $\RInc_k(2\times n)\setminus \Inc_k(2 \times n)$ and increasing tableaux in $\Inc_{k-1}(2 \times n)$. While this bijection does not preserve the descent set, by considering the change of the descent sets, we get a recurrence formula for $R_q(n,k)$ in terms of $S_q(n,k)$, and hence prove Theorem \ref{theo_rqnk} by applying Theorem \ref{theo_sqnk}. In Section 3 we give a bijection $ \Phi: \RInc_k(2\times n) \mapsto \RInc_k(2\times n)$ such that for any $T\in \RInc_k(2\times n),$
$$\maj(\Phi(T))=\amaj(T)+n-k.$$
Thus we proved Theorem \ref{theo_r2qnk} by showing that
$$\sum_{T \in \RInc_k(2 \times n)}q^{\maj(T)}=q^{n-k}\cdot\sum_{T \in \RInc_k(2 \times n)}q^{\amaj(T)}.$$
Finally in Section 4 we review some related work on enumerating major index of Schr\"{o}der paths.

\section{Counting major index for $\RInc_k(2 \times n)$}


\begin{theo} \label{theo_rnk}For any positive integer $n,k$, we have
\begin{eqnarray}
r(n,k)=s(n,k)+s(n,k-1).\label{recur_rnk}
\end{eqnarray}
\end{theo}

\pf Note that $\RInc_k(2 \times n)$ can be split into two disjoint sets: those that are increasing tableaux ($\Inc_k(2 \times n)$) and those that contains at least one column of identical entries ($\RInc_k(2\times n)\setminus \Inc_k(2 \times n)$). By definition the former one is counted by $s(n,k)$. And we will prove \eqref{recur_rnk} by providing a bijection between $\RInc_k(2\times n)\setminus \Inc_k(2\times n)$ and $\Inc_{k-1}(2\times n)$.

For any SSYT $T$, we use $T_{i,j}$ to denote the entry in row $i$ and column $j$ of $T$. Now we define $f$: $\RInc_k(2\times n)\setminus \Inc_k(2\times n)\xrightarrow{} \Inc_{k-1}(2\times n)$ as follows.
Given $T\in \RInc_k(2\times n)\setminus \Inc_k(2\times n)$, find the minimal integer $j$ such that $T_{1,j}=T_{2,j}$, i.e., the $j$-th column is the leftmost column of $T$ with two identical entries. Now we first delete the entry $T_{2,j}$, then move all the entries on the right of $T_{2,j}$ one box to the left and set the last entry as $2n-k+1$, and define the resulting tableau to be $f(T)$. Note that for any $i, 1 \leq i \leq n$, $T_{1,i} \leq T_{2,i} <T_{2,i+1}$, therefore in $f(T)$, we have $f(T)_{1,i} < f(T)_{2,i}$, and there are only $k-1$ number appear twice in $f(T)$. Hence $f(T) \in \Inc_{k-1}(2\times n)$.

The map $f$ is reversible. Given $S\in \Inc_{k-1}(2 \times n)$. Find the rightmost column $j^{\prime}$ such that $S_{1,j^{\prime}+1}=S_{2,j^{\prime}}+1$. (If such a column does not exist, then set $j^{\prime}=0$.) Now we first delete the entry $S_{2,n}$, then move all the entries $S_{2,j^{\prime}+1}, S_{2,j^{\prime}+2}, \ldots, S_{2,n-1}$ one box to the right, and set $S_{2,j^{\prime}+1}=S_{1,j^{\prime}+1}$. We denote the resulting tableau as $T$. It is obvious that $T=f^{-1}(S)\in \RInc_k(2\times n)\setminus \Inc_k(2\times n).$\qed

Figure \ref{fig_bij_Inck-1} shows an example of $f$ with $T \in \RInc_3(2\times 5)\setminus \Inc_3(2\times 5)$ and $f(T) \in \Inc_2(2\times 5)$.

\begin{figure}[ht]
\centering
$T:$\ \ \ytableaushort{13456,23467} $\ \ \ \ \ \ \ \; \mapsto\ \ \ \ \ \ f(T):$\ \  \ytableaushort{13456,24678}
\caption{An example of $f$ with $T \in \RInc_3(2\times 5)\setminus \Inc_3(2\times 5)$ and $f(T) \in \Inc_2(2\times 5)$.}\label{fig_bij_Inck-1}
\end{figure}

\begin{theo} \label{theo_recur_rqnk}
For any positive integer $n,k$ with $k<n$, we have
\begin{equation}
R_q(n,k)=S_q(n,k)+S_q(n,k-1)+(1-q^{2n-k})(S_q(n-1,k-1)+S_q(n-1,k-2)).\label{eq_recur_rqnk}
\end{equation}
\end{theo}
\pf Given $T \in \RInc_k(2\times n)$. We have either $T \in \Inc_k(2\times n)$ or $T \in \RInc_k(2\times n)\setminus \Inc_k(2\times n)$. The sum of $q^{\maj(T)}$ over all SSYTs in $\Inc_k(2\times n)$ is exactly $S_q(n,k)$. For all SSYTs in $\RInc_k(2\times n)\setminus \Inc_k(2\times n)$, there are two cases.
\begin{itemize}
\item[1.] If $T_{1,n}=T_{2,n}$. In this case we have $T_{1,n}=T_{2,n}=2n-k$ and $2n-k \notin D(T)$. We will show that the sum of $q^{\maj(T)}$ over all these tableaux is $S_q(n-1,k-1)+S_q(n-1,k-2)$.
\begin{itemize}
\item[1)] The $n$-th column is the only column of $T$ with identical entries. In this case the last column of $T$ consist of two identical entries $2n-k$ and $2n-k \notin D(T)$. And the sum of $q^{\maj(T)}$ over these tableaux is $S_q(n-1,k-1)$.

\item[2)] There is at least one column with identical entries in $T$ besides the $n$-th column. Now let $T^\prime$ be the tableau obtained by deleting the last column from $T$.  Clearly we have $T' \in \RInc_{k-1}(2 \times (n-1)) \setminus \Inc_{k-1}(2 \times (n-1))$.
There are two cases for the last column of $T'$.
\begin{itemize}
\item[a)] If $T'_{1,n-1}\neq T'_{2,n-1}$, then $f(T') \in \Inc_{k-2}(2 \times (n-1))$ with $f(T’)_{1,n-1}<2n-k-1$, and
$\maj(T) = \maj(T')=\maj(f(T'))$;
    
\item[b)] If $T'_{1,n-1}= T'_{2,n-1}$, we have $T_{1,n-1}= T_{2,n-1}=2n-k-1$. Since $T_{1,n}= T_{2,n}=2n-k$, we have $2n-k-1 \in D(T)$ but $2n-k-1 \notin D(T')$, and all the other descents of $T'$ are also descents of $T$. Thus
    \begin{equation}\label{eq1}
    \maj(T')=\maj(T)-(2n-k-1) .
    \end{equation}
    Moreover when we apply $f$ to $T'$ we have $2n-k-1 \in D(f(T'))$ but $2n-k-1 \notin D(T')$, and all the other descents of $T'$ are also descents of $f(T')$. Therefore we have
    \begin{equation}\label{eq2}
    \maj(f(T'))=\maj(T')+(2n-k-1).
    \end{equation}
Combining \eqref{eq1} and \eqref{eq2} we know that $f(T') \in \Inc_{k-2}(2 \times (n-1))$ with $f(T’)_{1,n-1}=2n-k-1$, and $\maj(f(T'))=\maj(T)$.
\end{itemize}
Thus the sum of  $q^{\maj(T)}$ over these tableaux of case a) and b) is $S_q(n-1,k-2)$.
\end{itemize}

\item[2.] If $T_{1,n} \ne T_{2,n}$. In this case applying $f$ to $T$ obtains a tableau in $\Inc_k(2
\times n)$ in which $T_{1,n} \ne 2n-k$.  The sum of $q^{\maj(T)}$ over all these tableau is $S_q(n,k)$ less $q^{2n-k}$ times the sum of
$q^{\maj(T)}$ over tableau in $\Inc_{k-1}(2 \times n)$ with $T_{1,n} = 2n-k$.
If we remove the last column of such a tableau (which is a reversible
operation since $T_{2,n} = 2n-k+1$) we obtain a tableau in either
$\Inc_{k-2}(2 \times (n-1))$ if $T_{2,n-1} = 2n-k$, or in $\Inc_{k-1}(2 \times
(n-1))$ if $T_{2,n-1} = 2n-k-1$. Therefore we have that the sum of $q^{\maj(T)}$ over these talbeaux is $S_q(n,k-1)-q^{2n-k}(S_q(n-1,k-1)+S_q(n-1,k-2))$.
\end{itemize}

Combining Case 1) and Case 2) we have
$$\sum_{T\in \RInc_k(2\times n)\setminus \Inc_k(2\times n)}q^{\maj(T)}=S_q(n,k-1)+(1-q^{2n-k})(S_q(n-1,k-1)+S_q(n-1,k-2)).$$
Hence \eqref{eq_recur_rqnk} is proved. \qed

\noindent  \emph{Proof of Theorem \ref{theo_rqnk}:}
From Theorem \ref{theo_sqnk} we have
\begin{equation*}
S_q(n,k)+S_q(n,k-1)=\frac{q^{n+k(k-1)/2}}{[n-k+1]}{2n-k \brack k}{2n-2k \brack n-k}.
\end{equation*}
Applying the above equation to Theorem \ref{theo_recur_rqnk} we have
\begin{eqnarray*}
R_q(n,k)
&=&\frac{q^{n+k(k-1)/2}}{[n-k+1]}{2n-k \brack k}{2n-2k \brack n-k}
   +(1-q^{2n-k})\frac{q^{n-1+(k-1)(k-2)/2}}{[n-k+1]}{2n-k-1 \brack k-1}{2n-2k \brack n-k}\\
&=&\frac{q^{n+k(k-3)/2}}{[n-k+1]}{2n-k \brack k}{2n-2k \brack n-k}.
\end{eqnarray*}\qed

\section{Counting amajor index for $\RInc_k(2 \times n)$}

In this section we will prove Theorem \ref{theo_r2qnk} by showing that

\[\sum_{T \in \RInc_k(2 \times n)}q^{\maj(T)}=q^{n-k}\cdot\sum_{T \in \RInc_k(2 \times n)}q^{\amaj(T)}.\]

Our main idea is to establish a bijection $\Phi: \RInc_k(2\times n) \mapsto \RInc_k(2\times n)$ which satisfies
$$\maj(\Phi(T))=\amaj(T)+n-k.$$

Before establishing the map we need some definitions. We say a row-increasing tableau $T$ is \emph{prime} if for each integer $j$ satisfies $T_{1,j+1}=T_{2,j}+1$, $T_{2,j+1}$ also appears in row 1 in  $T$. We use \emph{$\pRInc^m_k(\lambda)$} to denote the set of all prime row-increasing tableaux in $\RInc^m_k(\lambda)$.

For each $T\in \pRInc^m_k(2 \times n)$, we define two $k$-element sets $A$ and $B$ as the following: $A=\{a_1, a_2, \ldots, a_k\}_\leq$ is the set of numbers that appear twice in $T$. $B=\{b_1, b_2, \ldots, b_k\}$, here $b_i$ is the the number appears immediately left of $a_i$ in the second row of $T$ in cyclic order (if $a_1=T_{2,1}$, then $b_k=T_{2,n}$($=m+2n-k$)). Let $g(T)$ be the tableau of shape $2 \times n$ obtained by first deleting all elements in $A$ from the first row and then inserting all elements  in $B$ into the first row and list them in increasing order, and keep the entries in row 2 unchanged. (See Figure \ref{fig_g1} and Figure \ref{fig_tilde(T)} for examples.)

\begin{lemma}\label{lem_g}
The map $g$ is an injection from $\pRInc^m_k(2\times n)$ to $\RInc^m_k(2\times n)$ which satisfies the following:
\begin{itemize}
\item[1)] If $T_{2,1}$ appears only once in $T$, then $g(T)_{1,i+1} \leq g(T)_{2, i}$ for each $i, 1 \leq i \leq n-1$;
\item[2)] $T_{2,1}$ appears twice in $T$ if and only if $g(T)_{1,n}=g(T)_{2,n}$.
\end{itemize}
\end{lemma}
\pf  It is obvious that the map $g$ is invertible. Next we will prove that for any $T \in \pRInc^m_k(2\times n)$, $g(T)\in \RInc^m_k(2\times n)$ and satisfies corresponding conditions according to the following cases.

\begin{itemize}
\item[1.] When $T_{2,1}$ appears only once in $T$.

In this case for each $i, 1 \leq i \leq n$, $b_i$ is immediately to the left of $a_i$ in the second row of $T$, hence $b_i <a_i$.  Note that it is impossible that $T_{1,i+1} > T_{2, i}+1$ (in which case there will be no place for the number $T_{2, i}+1$), we now first prove that $g(T)_{1,i+1} \leq g(T)_{2, i}$ for each $i, 1 \leq i \leq n-1$ according to the following two cases.

\begin{itemize}
\item[a)] If $T_{1,i+1} \leq T_{2, i}$, we have $g(T)_{1,i+1} \leq T_{1,i+1} \leq T_{2,i} =g(T)_{2,i}$, for each $i, 1 \leq i \leq n-1$.

\item[b)] If $T_{1,i+1} = T_{2, i}+1$, according to the definition of prime row-increasing tableaux we know that in this case $T_{2,i+1}$ appears twice in $T$. Suppose there are exactly $x$ numbers among $T_{2,2}, T_{2,3}, \ldots, T_{2,i}$ that appear twice in $T$ ($x$ could be 0), then under the map $g$, exactly $x$ numbers are deleted from the left of $T_{1,i+1}$  in the first row and then $x+1$ numbers (including $T_{2,i}$) are inserted to the left of $T_{1,i+1}$, therefore we have $g(T)_{1,i+1}=T_{2,i}=g(T)_{2,i}$ for each $i, 1 \leq i \leq n-1$.
\end{itemize}

Moreover, since $g(T)_{1,i}<g(T)_{1,i+1}$ for each $i, 1 \leq i \leq n-1$, we know that $g(T)$ is strictly increasing in each column. And it is easy to check that $T_{1,1}=g(T)_{1,1}$ since $T$ is prime. Therefore we have $g(T) \in \RInc^m_k(2\times n)$ and $g(T)_{1,i+1}  \leq g(T)_{2,i}$ for each $i, 1 \leq i \leq n-1$. (See Figure \ref{fig_g1} for an example.)

\begin{figure}[h]
\centering
$T:$\ \ \ytableaushort{57{\textbf{8}}{10}{11}{\textbf{12}},
{6}{\textbf{8}}{9}{\textbf{12}}{13}{14}} $\ \  \ \ \; \underrightarrow{g}\ \ \ \
g(T):$\ \
\ytableaushort{{5}{\textbf{6}}{\textbf{7}}{\textbf{9}}{{10}}{{11}},
{\textbf{6}}{8}{\textbf{9}}{12}{13}{14}}
\label{fig_g1}
\caption{An example of the map $g$ with $T \in \pRInc_2^{4}(2\times 6)$ and $T_{2,1}$ appears only once.}
\end{figure}

\item[2.] When $T_{2,1}$ appears twice in $T$.

In this case we have $a_1=T_{2,1}$ and $b_1=T_{2,n}$. Let $\tilde{T}$ be the tableau of shape $2 \times n$ obtained from $T$ by first deleting $\{a_2, \ldots, a_k\}$ from the first row and then inserting $\{b_2, \ldots, b_k\}$ into the first row and list them in increasing order, and keep the entries in row 2 unchanged. (See Figure \ref{fig_tilde(T)} for an example.)

\begin{figure}[h]
\centering
$T:$\ \ \ytableaushort{1{\textbf{2}}{4}{5}{\textbf{6}}{\textbf{9}},{\textbf{2}}{{3}}{\textbf{6}}{7}{{8}}{\textbf{9}}} $\ \  \ \ \;\underrightarrow{g}\ \ \ \  g(T):$\ \ \ytableaushort{{{1}}{\textbf{3}}{4}{{5}}{\textbf{8}}{\textbf{9}},{2}{\textbf{3}}{6}{7}{\textbf{8}}{\textbf{9}}}\\ \vskip3mm
%
$\tilde{T}:$\ \ \ytableaushort{{1}{\textbf{2}}{\textbf{3}}{4}{{5}}{\textbf{8}},{\textbf{2}}{\textbf{3}}{6}{{7}}{\textbf{8}}{9}}
\caption{An example of the map $g$ with $T \in \pRInc_2^{0}(2\times 6)$ and $T_{2,1}$ appears twice.}
\label{fig_tilde(T)}
\end{figure}

Suppose $T_{2,1}=T_{1,j}, 1 \leq j \leq n$, then we have $g(T)_{1,i}=\tilde{T}_{1,i}$ for each $i, 1 \leq i \leq j-1$；
$g(T)_{1,i}=\tilde{T}_{1,i+1}$ for each $i, j \leq i \leq n-1$； and $g(T)_{1,n}={T}_{2,n}=g(T)_{2,n}$. Similar to the argument in case 1, we can prove that $\tilde{T}_{1,i+1} \leq \tilde{T}_{2,i}=T_{2,i}$ for each $i, 1 \leq i \leq n-1$.
Hence we have that $g(T)_{1,i} \leq g(T)_{2,i}$ for each $i, 1 \leq i \leq n-1$, i.e., $g(T)$ is weakly increasing in each column. And it is obvious from the definition of $g$ that $g(T)_{1,n}=g(T)_{2,n}=T_{2,n}$ if and only if $T_{2,1}$ appears twice in $T$.
\end{itemize}\qed

\begin{prop}\label{prop_maj=amaj+mnk}
For each $T \in \pRInc^m_k(2\times n)$  we have
\begin{equation}\label{eq_maj=amaj+mnk}
\maj(g(T))=\left\{
\begin{aligned}
&\amaj(T)+n-k,\ \ \ \ \ \text{if} \ T_{1,1}=T_{2,1};\\
&\amaj(T)+m+n-k,\ \ \ \ \ \text{if} \  T_{1,1}\ne T_{2,1}.
\end{aligned}
\right.
\end{equation}
\end{prop}
\pf Given $T \in \pRInc^m_k(2\times n)$, let $T^0$ be the skew shape tableau obtained by deleting the numbers in $A$ from row 1 of $T$ and “push” all the remaining numbers to the right, and keep the second row unchanged.
(See Example \ref{Eg_T0} for an example.) We will prove Equation \eqref{eq_maj=amaj+mnk} by verifying the following facts.

\begin{itemize}
\item[1.] $D(g(T))\setminus D(T^0)=A(T)\setminus A(T^0),$  and therefore
\[\maj(g(T))-\maj(T^0)=\amaj(T)-\amaj(T^0);\]

For each $i, 1 \leq i \leq n$, $i \in A(T)\setminus A(T^0)$ if and only if that in $T$ $i$ appears in row 2 and $i+1$ appears in both rows. In this case $i$ is immediately to the left of $i+1$ in row 2 in $T$ and is inserted to the first row of $T^0$ to get $g(T)$, hence $i \in D(g(T))\setminus D(T^0)$.

On the other hand, for each $i, 1 \leq i \leq n$, $i \in D(g(T))\setminus D(T^0)$ if only if that in $g(T)$, $i+1$ appears in row 2 and $i$ appears in both rows. In this case $i+1$ is immediately to the right of $i$ in $g(T)$, and hence $i+1$ appears in both rows of $T$ and $i$ appears in row 2 of $T$, therefore we have $i \in A(T)\setminus A(T^0)$.

\item[2.]
We first prove the case when $T_{1,1}\ne T_{2,1}$, $\maj(T^0)=\amaj(T^0)+m+n-k.$

Suppose $|D(T^0)|=d$ for some positive integer $d$, then we have $|A(T^0)|=d-1$. Moreover, suppose the descents of $T^0$ appears in columns $k+x_1, k+x_2,$ $\ldots,$ $k+x_d$ in row 1, and the ascents of $T^0$ appears in columns $y_1, y_2,$ $\ldots,$ $y_{d-1}$ in row 2. Here $d,$ $x_1, \ldots, x_d$, $y_1, \ldots y_{d-1}$ are all positive integers and $x_d=n-k$. It is not hard to check that $T^0$ is uniquely determined by the two sets $X=\{x_1, x_2, \ldots, x_d\}_\leq$ and $Y=\{y_1, y_2, \ldots, y_{d-1}\}_\leq$. And we have
\begin{eqnarray*}
D(T^0)&=&\{m+x_1, m+x_2+y_1, \ldots, m+x_{d}+y_{d-1}\};\\
A(T^0)&=&\{m+x_1+y_1, m+x_2+y_2, \ldots, m+x_{d-1}+y_{d-1}\}.
\end{eqnarray*}
Therefore we have
\[
\maj(T^0)-\amaj(T^0)
=m+x_d=m+n-k.
\]
Similarly we can prove that when $T_{1,1}= T_{2,1}$, $\maj(T^0)=\amaj(T^0)+n-k.$
\end{itemize}
Combining the above two facts we get Equation \eqref{eq_maj=amaj+mnk}.
\qed

\begin{exam}\label{Eg_T0}
Figure \ref{fig_amj-maj-mnk} shows an example of $T$, $T^0$ and $g(T)$ with $n=6$, $m=4$, and $k=2$. Here we have $A(T)=\{7,8,12\}$, $A(T^0)=\{8\}$, $D(T^0)=\{6,10\}$, $D(g(T))=\{6,7,10,12\}$, $d=2$, $X=\{2,4\}$ and $Y=\{2\}$.

\begin{figure}[h]
\centering
$T:$\ \ \ytableaushort{56{\textbf{8}}{9}{10}{\textbf{13}},
{\color{red}{7}}{\color{red}\textbf{8}}{11}{\color{red}{12}}{\textbf{13}}{14}} $\ \  \ \ \; \underrightarrow{g}\ \ \ \
g(T):$\ \
\ytableaushort{{5}{\color{blue}{6}}{\color{blue}\textbf{7}}{9}{\color{blue}{10}}{\color{blue}\textbf{{12}}},
{\textbf{7}}8{11}{\textbf{12}}{13}{14}}
\vskip 5mm

$T^0:$\ \ \ytableaushort{{}{}{5}{\color{blue}{6}}{9}{\color{blue}{10}},{7}{\color{red}8}{11}{12}{13}{14}}
\caption{An example of the map $g$ with $T \in \pRInc_2^{4}(2\times 6).$}
\label{fig_amj-maj-mnk}
\end{figure}

\end{exam}
Now we are ready to prove the general case.

\begin{theo}\label{mainthm_Phi}
There is a bijection  $\Phi: \RInc_k(2\times n) \vdash \RInc_k(2\times n)$ such that for any $T \in \RInc_k(2\times n)$, the second row of $\Phi(T)$ is identical with the second row of $T$, and
\begin{equation}\label{eq_amaj=maj+n-k}
\maj(\Phi(T))=\amaj(T)+n-k.
\end{equation}
\end{theo}
\pf
Given $T \in \RInc_k(2\times n)$, there is a unique way to decompose $T$ into prime row-increasing tableaux: suppose $i_1, i_2, \ldots, i_{l-1}$ are all the positive integers such that $T_{2, i_j}+1=T_{1, i_j+1}$ and $T_{2, i_j+1}$ appear only once in $T$, we break $T$ between column $i_j$ and $i_j+1$ and get a decomposition $T_1 T_2 \cdots T_l$ of $T$ into prime row-increasing tableaux. Now set $\Phi(T)=g(T_1)g(T_2)\cdots g(T_l)$ (see Example \ref{Eg_general}).

From Lemma \ref{lem_g} we know that $\Phi(T) \in \RInc_k(2\times n)$. Next we show that Equation \ref{eq_amaj=maj+n-k} holds. Suppose $T_j \in \pRInc^{m_j}_{k_j}(2\times n_j)$ for integers $m_j, k_j, n_j$ with $m_j, k_j \geq 0$, and  $n_j>0$. Then we have
$$n_1+n_2+\cdots+n_l=n, \ \ k_1+k_2+\cdots+k_l=k.$$
The smallest entry of $T_i$ is $m_i +1$ with $m_1=0$, and
$$m_j=2(n_1+n_2+\cdots+n_{j-1})-(k_1+k_2+\cdots+k_{j-1}), \ \ 2 \leq j \leq l.$$
And the largest entry of $T_j$ is  $m_{j+1}$ for each $j, 1 \leq j \leq l-1$,

It is easy to check that
$$A(T)=A(T_1)\cup A(T_2)\cup \cdots \cup A(T_l)\cup \{m_2, m_3, \ldots, m_l\}$$
and $$\amaj(T)=
\sum_{j=1}^{l} \amaj(T_j)+\sum_{j=2}^{l}m_j.$$

Moreover we have
$$D(\Phi(T))=D(g(T_1))\cup D(g(T_2))\cup \cdots \cup D(g(T_l)).$$

From the definition of prime row-increasing tableaux we know it is impossible that the two numbers of the first column of $T_j$ is identical when $j >1$. Moreover, since $T \in \RInc_k(2\times n)$ we have $m_1=0$. Hence from Proposition \ref{prop_maj=amaj+mnk} we know that for each $j, 1 \leq j \leq l,$
$$\maj(g(T_j))=\amaj(T_j)+m_j+n_j-k_j,$$
always holds. Therefore we have
\begin{eqnarray*}
\maj(\Phi(T))&=&\sum_{j=1}^l{\maj(g(T_j))}
=\sum_{j=1}^l(\amaj(T_j)+m_j+n_j-k_j)\\
&=&\sum_{j=1}^l\amaj(T_j)+\sum_{j=1}^l m_j+n-k \\
&=&\amaj(T)+n-k.
\end{eqnarray*}

It remains to show that $\Phi$ is a bijection. Since $g$ is an injection from $\pRInc^m_k(2\times n)$ to $\RInc^m_k(2\times n)$, it is sufficient to show that  given $S=\Phi(T) \in \RInc_k(2\times n)$, we can decompose $S$ at the right places to get $S_1 S_2 \cdots S_l$ such that $S_i=g(T_i)$ for each $i, 1 \leq i \leq l$. And the decomposition is as follows, we first find the largest number $j$ such that $S_{1,j}=S_{2,j}$, i.e., column $j$ is the rightmost column with identical entries. If such a column does not exist, we set $j=0$. Now for each $i, j\leq i \leq n-1$, we break $S$ between column $i$ and $i+1$ if $S_{1,i+1}>S_{2,i}$ (Note that when $0<j<n$,  $S_{1,j+1}=S_{2,j}+1$ always holds) and get a decomposition $S_1 S_2 \cdots S_t$ of $S$.  From Lemma \ref{lem_g} it is clear that such a decomposition guarantees $t=l$ and $S_i=g(T_i)$ for each $i, 1 \leq i \leq l$. Therefore we proved that $\Phi$ is a bijection.
\qed

\begin{exam}\label{Eg_general}
Let $n=13$, $k=6$, and $l=3$. Figure \ref{fig_general} shows an example of the map $\Phi$. Here we have $A(T)=\{3, 8, 9, 11, 13, 15, 17, 19\}$, $A(T_1^0)=\{3\}$,  $A(T_2^0)=\{11,13\}$, $A(T_3^0)=\emptyset$, $D(T_1^0)=\{1,5\}$,  $D(T_2^0)=\{10,12,14\}$, $D(T_3^0)=\{18\}$. $D(\Phi(T))=\{1, 5, 8, 10, 12, 14, 15, 18, 19\}$. Here we have $\amaj(T)=95$ and $\maj(\Phi(T))=102$.
\begin{figure}[h]
\centering
\ \ \ \ $T:$\ \ \ytableaushort{1{\textbf{2}}45{\textbf{6}}{\textbf{9}}{10}{12}{\textbf{13}}{14}{\textbf{16}}{18}{\textbf{20}},{{\textbf{2}}}{\color{red}{3}}{{\textbf{6}}}{7}{\color{red}{8}}{\color{red}{\textbf{9}}} {\color{red}{11}}{\color{red}{\textbf{13}}}{\color{red}{15}}{\textbf{16}}{\color{red}{17}}{\color{red}{19}}{\textbf{20}}}

\vskip 5mm

\ \ \ \ \ \ \ \ \ \ \ytableaushort{{}{}{}{\color{blue}{1}}4{\color{blue}{5}}{}{}{\color{blue}{10}}{\color{blue}{12}}{\color{blue}{14}}{}{\color{blue}{18}},{2}{\color{red}{3}}6789{\color{red}{11}}{\color{red}{13}}{15}{16}{17}{19}{20}}
\vskip 5mm

$\Phi(T):$\ \ \ytableaushort{{\color{blue}{1}}{\textbf{3}}{4}{\color{blue}{5}}{\color{blue}{\textbf{8}}}{\textbf{9}}{\color{blue}{10}}{\textbf{11}}{\color{blue}{12}}{\color{blue}{14}}{\color{blue}{\textbf{15}}}{\color{blue}{18}}{\color{blue}{\textbf{19}}},2{\textbf{3}}67{\textbf{8}}{\textbf{9}}{\textbf{11}}{13}{\textbf{15}}{16}{17}{\textbf{19}}{20}}
\caption{An example of the map $\Phi$ with $n=13$, $k=6$, and $l=3$.}\label{fig_general}
\end{figure}
\end{exam}

\section{Major index of Schr\"{o}der $n$-paths}\label{sec2.3}

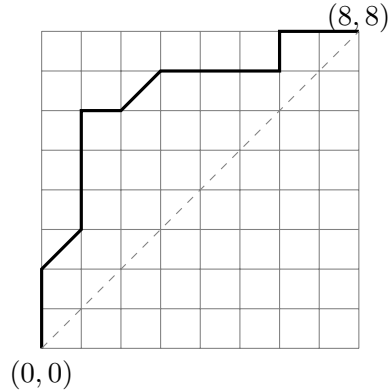
\begin{figure}[ht]
\centering
\begin{tikzpicture}[thick,domain=0:8,mark=]
\draw[ thin,gray,step=15pt] (0,0) grid (120pt,120pt);
\draw[yshift=0cm,very thick] plot coordinates %
{(0pt,0pt) (0pt,30pt) (15pt,45pt) (15pt,90pt) (30pt,90pt) (45pt,105pt) (90pt,105pt)(90pt,120pt) (120pt,120pt)};
\draw [dashed,thin,gray](0,0)--(120pt,120pt);
\node at (0pt,-10pt) [rotate=0] {$(0,0)$};
\node at (120pt,125pt) [rotate=0] {$(8,8)$};
\end{tikzpicture}
\caption{A Schr\"{o}der $P$ with $\omega(P)=00100021222022$.}\label{Schroder}
\end{figure}

Let $P$ be a Schr\"{o}der $n$-path that goes from the origin $(0,0)$ to $(n,n)$ with $k$ $F$ steps, we can associate with $P$ a word $w=w(P)=w_1 w_2 \cdots w_{2n-k}$ over the alphabet $\{0,1,2\}$ with exactly $k$ $1$'s. (See Firgure \ref{Schroder}.)

We say that $w$ has a \emph{descent} in position $i, 1 \leq i \leq n-1$, if $w_i > w_{i+1}$. The \emph{descent set} of $w$ is the set of all positions of the descents of $w$, $D(w)=\{i: 1 \leq i \leq n-1, w_i > w_{i+1}\}$. The \emph{major index} of $w$ is defined as $\maj(w)=\sum_{i \in D(w)}i$. And define $\maj(P)=\maj(w(P))$.

In \cite{Bonin_Shapiro_Simion}, Bonin, Shapiro and Simion study the major index for Schr\"{o}der paths and gave the following result:
\begin{equation}\label{maj_Sch}
\sum q^{\maj(P)}=\frac{1}{[n-k+1]}{2n-k \brack k}{2n-2k \brack n-k},
\end{equation}
here the sum is over all Schr\"{o}der $n$-paths with exactly $k$ $F$ steps.

Note that Equation \eqref{maj_Sch} differs from Equation \eqref{rqnk} and \eqref{r2qnk} only by a factor of a power of $q$. Readers might wonder if there is some simple explanation on these relations. In fact there is an obvious bijection $\theta$ between SSYTs in $\RInc_k(2 \times n)$ and Schr\"{o}der words of length $2n-k$ that contain exactly $k$ 1’s. Given $T \in \RInc_k(2 \times n)$, we read the numbers from $1$ to $2n-k$ in increasing order, if $i$ appears only in row 1 (2), we set $w_i=0 (2)$, otherwise we set $w_i=1$, and define $\theta(T)=w_1 w_2 \cdots w_{2n-k}$. A naive thinking is that if $i$ is an ascent of $T$, then $i$ is a descent of $\theta(T)$, but this is not always true. When $i$ and $i+1$ both appear in row 1 and row 2 of $T$, then we have $i \in A(T)$ but $i \notin D(\theta(T))$. It would be interesting if one can find a simple combinatorial explanation on the relations between descent (ascent) sets of a row-increasing tableaux in $\RInc_k(2 \times n)$ and the descent sets of the correspongding Schr\"{o}der $n$-paths with $k$ $F$ steps.

\vskip 2mm \noindent{\bf Acknowledgments.} We thank the referee for a very careful reading and many important suggestions adopted in this revised version. This work is partially supported by National Natural Science Foundation of China (No. 11871223) and the Science and Technology Commission of Shanghai Municipality (No. 18ZR1411700 and No. 18dz2271000).

\end{document}